\documentclass[12pt,reqno]{amsart}

\usepackage{amssymb}

\usepackage{eucal}

\input{diagrams}

\font\sss=cmss8

\def\cG{{\mathcal G}}

\def\BZ{{\mathbb Z}}

\def\sD{\mbox{\sf D}}
\def\sE{\mbox{\sf E}}

\def\sK{\mbox{\sf K}}

\def\Ab{\mbox{\sf Ab}}

\def\ast{{\textstyle *}}

\def\Coker{\operatorname{Coker}}
\def\D{\sD}

\def\Ext{\operatorname{Ext}}

\def\GExt{\operatorname{Ext}_{\cG}}

\def\H{\operatorname{H}}

\def\Hom{\operatorname{Hom}}
\def\id{\operatorname{id}}
\def\Image{\operatorname{Im}}

\def\Ker{\operatorname{Ker}}

\def\LTensor{\stackrel{\operatorname{L}}{\otimes}}
\def\Mod{\mbox{\sf Mod}}

\def\opp{\operatorname{op}}

\def\Proj{\mbox{\sf Pro}}

\def\RHom{\operatorname{RHom}}

\def\TExt{\widehat{\operatorname{Ext}}}

\numberwithin{equation}{part}

\newtheorem{Lemma}{Lemma}[section]
\newtheorem{Theorem}[Lemma]{Theorem}
\newtheorem{Proposition}[Lemma]{Proposition}

\theoremstyle{definition}
\newtheorem{Definition}[Lemma]{Definition}
\newtheorem{Setup}[Lemma]{Setup}

\newtheorem{Construction}[Lemma]{Construction}
\newtheorem{Remark}[Lemma]{Remark}

\def\R{A}

\def\KProjR{\sK(\Proj\,\R)}
\def\KProjRsmall{\mbox{\sss K}(\mbox{\sss Pro}\,\R)}
\def\EnochsR{\sE(\R)}
\def\EnochsRsmall{\mbox{\sss E}(\R)}

\def\inc{e_{\ast}}
\def\adj{e^!}
\def\res{\operatorname{res}}

\def\Rlm{$\R$-left-module}

\def\M{M}          
\def\N{N}       

\def\Gorproj{G}         
\def\Gorprojtil{\widetilde{G}}
\def\projmod{Q}         
\def\flatmod{F}         

\def\kercomplx{K}       
\def\kercomplxprime{K^{\prime}}
\def\projcomplx{P}      
\def\projcomplxtil{\widetilde{P}}

\def\Enochscomplx{E}    

\def\Enochscomplxtil{F} 
\def\Enochscomplxtilprime{F^{\prime}}
\def\Tatecomplx{T}      

\def\Enochscomplxmap{e} 
 
\def\Gorprojmap{g}   
\def\Gorprojtilmap{\widetilde{g}}
\def\Tatecomplxmap{t}   

\def\Rm{$\R$-module}

\def\injmod{J}          
\def\injcomplx{I}       

\def\generalR{A}
\def\matlisgeneralR{B}

\def\generalRm{$\generalR$-module}
\def\generalRlm{$\generalR$-left-module}

\def\EnochsgeneralR{\sE(\generalR)}
\def\KProjgeneralR{\sK(\Proj\,\generalR)}

\begin{document}

\title[Gorenstein resolutions]
{Existence of Gorenstein projective resolutions}

\author{Peter J\o rgensen}
\address{Department of Pure Mathematics\\
         University of Leeds\\
         Leeds LS2 9JT\\
         United Kingdom}
\email{popjoerg@maths.leeds.ac.uk}
\urladdr{http://www.maths.leeds.ac.uk/\~{ }popjoerg}


\keywords{Dualizing complex, Bousfield localization, Gorenstein
homological algebra, Gorenstein projective precover, Tate cohomology}

\subjclass[2000]{13D02, 16E05, 18G25, 20J06}

\begin{abstract}

Gorenstein rings are important to mathematical areas as diverse as
algebraic geometry, where they encode information about singularities
of spaces, and homotopy theory, through the concept of model
categories. 

In consequence, the study of Gorenstein rings has led to the advent of
a whole branch of homological algebra, known as Gorenstein homological
algebra.

This paper solves one of the open problems of Gorenstein homological
algebra by showing that so-called Gorenstein projective resolutions
exist over quite general rings, thereby enabling the
de\-fi\-ni\-ti\-on of a Gorenstein version of derived functors.

An application is given to the theory of Tate cohomology.

\end{abstract}

\maketitle

\setcounter{section}{-1}
\section{Introduction}
\label{sec:introduction}

Gorenstein rings are important mathematical objects originating in the
work of Grothendieck and his pupils.  The study of Gorenstein rings
has given rise to a whole branch of homological algebra known as
Gorenstein homological algebra, to which this paper is a contribution.
The main point is an existence proof for Gorenstein projective
resolutions; this central item has been lacking from the theory until
now.  However, before going into details, let me start the exposition
on a classical note.

{\em Homological algebra} is one of the most versatile mathematical
machines ever invented.  It has an impact on most parts of
mathematics --- algebra, geometry, number theory...  One of the central
items of the theory are {\em projective resolutions}.  An augmented
projective resolution of a module $\M$ is an exact sequence
\[
  \cdots \longrightarrow \projcomplx_2
  \longrightarrow \projcomplx_1
  \longrightarrow \projcomplx_0
  \longrightarrow \M
  \rightarrow 0
\]
where the $\projcomplx_i$ are projective modules.  If each module $\M$
over a ring $\R$ has a projective resolution $\projcomplx$ with
$\projcomplx_i = 0$ for $i \gg 0$, then $\R$ is called {\em regular}.
For an algebraist, one of the points of homological algebra is that it
can be used to study non-regular rings in terms of their deviation
from regularity; this makes it possible to use the easier regular
case as a frame of reference for understanding more complicated rings.

{\em Gorenstein homological algebra} takes a parallel approach and
considers {\em Gorenstein projective resolutions}.  I will give the
precise definition in a moment, but the point is that if each module
$\M$ over a ring $\R$ has a Gorenstein projective resolution
$\Gorproj$ with $\Gorproj_i = 0$ for $i \gg 0$, then $\R$ is {\em
Gorenstein}, that is, it has finite injective dimension as a module
over itself.  Gorenstein homological algebra goes back to Auslander
and Bridger, see \cite{AuslanderBridger}, and substantial
contributions are due to Enochs and his coauthors, see
\cite{EJLR}, \cite{EJX}, \cite{ELR}, and several other papers.

Gorenstein homological algebra plays a role in algebraic geometry, see
\cite{EvansGriffith1}, \cite{EvansGriffith2}, \cite{Hartshorne};
commutative ring theory, see \cite{AuslanderBuchweitz},
\cite{AvrMart}, \cite{AvrFoxPLMS}, \cite{HHGorensteinDerived},
\cite{HHGorensteinHomDim}, \cite{Veliche}; homotopy theory, see
\cite{Hovey}; and number theory, see \cite{Venjakob}.


To describe {\em the contents of this paper}, let me make a foray into
the definitions of Gorenstein homological algebra.  Gorenstein
projective resolutions are defined in terms of {\em Gorenstein
projective modules}.  These are modules of the form $\Gorproj =
\Ker(\Enochscomplx^1 \longrightarrow \Enochscomplx^2)$ where
$\Enochscomplx$ is a {\em complete projective resolution}, that is, an
exact complex of projective modules which stays exact when one applies
the functor $\Hom(-,\projmod)$ for any projective module $\projmod$.
An augmented Gorenstein projective resolution of a module $\M$ is an
exact sequence
\begin{equation}
\label{equ:aug_Gor_proj_res}
  \cdots \longrightarrow \Gorproj_2
  \longrightarrow \Gorproj_1
  \longrightarrow \Gorproj_0
  \longrightarrow \M
  \rightarrow 0
\end{equation}
where the $\Gorproj_i$ are Gorenstein projective modules, which
stays exact when one applies the functor $\Hom(\Gorprojtil,-)$
for any Gorenstein projective module $\Gorprojtil$.  The complex
\begin{equation}
\label{equ:Gor_proj_res}
  \Gorproj = 
  \cdots \longrightarrow \Gorproj_2
  \longrightarrow \Gorproj_1
  \longrightarrow \Gorproj_0
  \longrightarrow 0
  \longrightarrow \cdots
\end{equation}
is then called a Gorenstein projective resolution of $\M$.

For the theory to be worth anything, it is a key question whether
Gorenstein projective resolutions really exist.  In other words, for a
given ring $\R$, it is important to determine whether each \Rm\ has a
Gorenstein projective resolution.  One could attempt to circumvent
this question by dropping the requirement that the complex
\eqref{equ:aug_Gor_proj_res} stay exact under the functor
$\Hom(\Gorprojtil,-)$.  This makes it easy to establish existence, and
such attempts have been made.  However, to do so misses an important
point: The purpose of requiring \eqref{equ:aug_Gor_proj_res} to stay
exact under the functor $\Hom(\Gorprojtil,-)$ is that this makes the
Gorenstein projective resolution \eqref{equ:Gor_proj_res} unique up to
chain homotopy, as one can easily check.  This in turn means that
\eqref{equ:Gor_proj_res} can be used to define the {\em Gorenstein
version of derived functors}.  Without the requirement that
\eqref{equ:aug_Gor_proj_res} stay exact under the functor
$\Hom(\Gorprojtil,-)$, any such definition fails, and the theory must
remain without derived functors; a one legged life.  Therefore, it is
a central question whether Gorenstein projective resolutions exist.

The corresponding questions of existence of so-called Gorenstein
injective and Gorenstein flat resolutions have recently been settled
in the affirmative in \cite{ELR} and \cite{EJLR}, but the Gorenstein
projective case has resisted the attacks of a number of authors
despite partial results in papers such as \cite{AuslanderBuchweitz},
\cite{EJX}, \cite{HHGorensteinHomDim}, and \cite{PJSpectra}.  The
state of the art up to now seems to be \cite[prop.\
2.18]{HHGorensteinHomDim}; this only gives that Gorenstein
projective resolutions exist over a Gorenstein ring.

However, the present paper proves the existence of Gorenstein
projective resolutions over much more general rings.  This is done by
showing that the resolutions exist under one simple assumption --- the
existence of a certain {\em adjoint functor} $\adj$ --- and by using
Bousfield localization to show that this assumption holds if the
ground ring has a {\em dualizing complex}.  This covers many rings
arising in practice.  For instance, any local ring of a scheme of
locally finite type over a field has a dualizing complex.  Other types
of rings are also covered; see remark \ref{rmk:dualizing}.

In fact, it may even be the case that the functor $\adj$ exists over
{\em any} ring and hence that Gorenstein projective resolutions exist
in general, but I do not know how to prove this.

After showing these results, I will give an application to the theory
of Tate cohomology.  This was originally defined for representations
of finite gro\-ups, but I will show, again under the assumption that
the adjoint functor $\adj$ exists, that it is possible to define Tate
$\Ext$ groups
\[
  \TExt^i(\M,\N)
\]
for any modules $\M$ and $\N$, so that classical Tate cohomology is
the special case $\TExt_{kG}^i(k,\N)$.  Moreover, Tate and
ordinary $\Ext$ groups will be shown to fit into a long exact sequence
\begin{equation}
\label{equ:long_exact_sequence}
 0 \rightarrow \GExt^1(\M,\N) 
 \longrightarrow \Ext^1(\M,\N)
 \longrightarrow \TExt^1(\M,\N) 
 \longrightarrow \cdots,
\end{equation}
where the $\GExt^i$ are Gorenstein $\Ext$ groups defined by
\[
  \GExt^i(\M,\N) = \H^i\!\Hom(\Gorproj,\N)
\]
where $\Gorproj$ is a Gorenstein projective resolution of $\M$.  The
$\GExt^i$ are precisely a Gorenstein version of derived functors.

A theory of Tate cohomology such as this was already accomplished in
\cite{AvrMart} and \cite{Veliche}, but only under the assumption that
$\M$ had a finite Gorenstein projective resolution, hence restricting
the real scope of the theory to Gorenstein rings.

The paper is organized as follows.  Section \ref{sec:adjoint} shows
the existence of the adjoint functor $\adj$ over rings with a
dualizing complex.  Section \ref{sec:resolution} shows the existence
of Gorenstein projective resolutions when $\adj$ exists.  And section
\ref{sec:Tate} defines Tate $\Ext$ groups, shows some simple properties,
and shows that the Tate $\Ext$ groups fit into the exact sequence
\eqref{equ:long_exact_sequence}.

\setcounter{section}{0}
\section{The adjoint functor $\adj$}
\label{sec:adjoint}

This section shows the existence of a certain adjoint functor $\adj$
over rings with a dualizing complex.

\begin{Remark}
\label{rmk:dualizing}
Dualizing complexes are popular gadgets in homological algebra.  I
shall give the precise definition in setup \ref{set:blanket1} for
noetherian commutative rings and in setup \ref{set:blanket1}' for
right-noetherian algebras over a field.  But I would like already here
to point out that many rings have dualizing complexes.

For instance, a noetherian local commutative ring has a dualizing
complex if and only if it is a quotient of a Gorenstein noetherian
local commutative ring, by the (deep) result \cite[thm.\
1.2]{Kawasaki}.  It follows that, as mentioned in the introduction,
any local ring of a scheme of locally finite type over a field has a
dualizing complex.  By the Cohen structure theorem, it also follows
that any complete noetherian local commutative ring does.

Some important types of non-commutative noetherian algebras are also
known to have dualizing complexes.  For example, complete semi-local
PI algebras do by \cite[cor.\ 0.2]{WuZhangDualizing}, and filtered
algebras do by \cite[cor.\ 6.9]{YekutieliZhang} if their associated
graded algebras are noetherian and connected, and either PI, FBN, or
with enough normal elements.
\end{Remark}

\begin{Definition}
\label{def:Enochs}
If $\generalR$ is a ring, then $\EnochsgeneralR$ denotes the class
of complete projective resolutions of \generalRlm s.  So a complex of
\generalRlm s $\Enochscomplx$ is in $\EnochsgeneralR$ if it consists of 
projective \generalRlm s, is exact, and has
$\Hom_{\generalR}(\Enochscomplx,\projmod)$ exact for each projective
\generalRlm\ $\projmod$.
\end{Definition}

\begin{Remark}
I will view $\EnochsgeneralR$ as a full subcategory of
$\KProjgeneralR$, the homotopy category of complexes of projective
\generalRlm s.  The inclusion functor will be denoted
\[
  \inc : \EnochsgeneralR \longrightarrow \KProjgeneralR.
\]
\end{Remark}

\begin{Setup}
\label{set:blanket1}
Let $\generalR$ be a noetherian commutative ring with a dualizing complex
$D$.  That is,
\begin{enumerate}

  \item  The cohomology of $D$ is bounded and finitely generated
         over $\generalR$.

  \item  The injective dimension $\id_{\generalR} D$ is finite.

  \item  The canonical morphism $\generalR \longrightarrow
         \RHom_{\generalR}(D,D)$ in the derived category $\D(\generalR)$ is an
         isomorphism. 

\end{enumerate}
\end{Setup}

\begin{Setup}
\label{set:I1}
Let $D \stackrel{\simeq}{\longrightarrow} \injcomplx$ be an injective
resolution so that $\injcomplx$ is a bounded complex.
\end{Setup}

See \cite[chp.\ V]{HartsResDual} for background on dualizing
complexes. 

\begin{Remark}
\label{rmk:flats_have_finite_pd}
Since $\generalR$ has a dualizing complex, it has finite Krull dimension
by \cite[cor.\ V.5.2]{HartsResDual}, so by \cite[Seconde partie, cor.\
(3.2.7)]{RaynaudGruson}, each flat \generalRm\ has finite projective
dimension.
\end{Remark}

The following lemma uses $\injcomplx$, the injective resolution of the
dualizing complex $D$.

\begin{Lemma}
\label{lem:technical}
Let $\projcomplx$ be a complex of projective \generalRm s.  Then
\begin{eqnarray*}
  & \mbox{$\Hom_{\generalR}(\projcomplx,\projmod)$ is exact for each
          projective \generalRm\ $\projmod$} & \\
  & \mbox{$\Leftrightarrow \injcomplx \otimes_{\generalR} \projcomplx$
          is exact.} &
\end{eqnarray*}
\end{Lemma}

\begin{proof}
\noindent
$\Rightarrow \;$  Suppose that $\Hom(\projcomplx,\projmod)$ is exact
for each projective module $\projmod$.  To see that $\injcomplx
\otimes \projcomplx$ is an exact complex, it is enough to see that
\[
  \Hom(\injcomplx \otimes \projcomplx,\injmod) 
  \cong \Hom(\projcomplx,\Hom(\injcomplx,\injmod))
\]
is exact for each injective module $\injmod$.  

It follows from \cite[thm.\ 1.2]{Lazard} that
$\Hom(\injcomplx,\injmod)$ is a bounded complex of flat modules.
Hence, $\Hom(\injcomplx,\injmod)$ is finitely built from flat modules
in the homotopy category of complexes of \generalRm s,
$\sK(\generalR)$, and so it is enough to see that
$\Hom(\projcomplx,\flatmod)$ is exact for each flat module $\flatmod$.

Since $\flatmod$ has finite projective dimension by remark
\ref{rmk:flats_have_finite_pd}, there is a projective resolution
$\projcomplxtil \stackrel{\simeq}{\longrightarrow} \flatmod$ with
$\projcomplxtil$ bounded.  Since $\projcomplx$ consists of projective
modules and both $\projcomplxtil$ and $\flatmod$ are bounded, this
induces a quasi-isomorphism
\[
  \Hom(\projcomplx,\projcomplxtil) \simeq \Hom(\projcomplx,\flatmod).
\]
So it is enough to see that $\Hom(\projcomplx,\projcomplxtil)$ is exact.  

But $\projcomplxtil$ is a bounded complex of projective modules, so it
is finitely built from projective modules, so it is enough to see that
$\Hom(\projcomplx,\projmod)$ is exact for each projective module
$\projmod$.  And this holds by assumption.

\smallskip

\noindent
$\Leftarrow \;$ Suppose that $\injcomplx \otimes \projcomplx$ is an
exact complex.  I must show that $\Hom(\projcomplx,\projmod)$ is exact
for each projective module $\projmod$.

First observe that by \cite[thm.\ (3.2)]{AvrFoxPLMS}, there is an
isomorphism
\[
  \projmod 
  \stackrel{\sim}{\longrightarrow} \RHom(D,D \LTensor \projmod).
\]
Of course, I can replace $D$ by $\injcomplx$ to get
\begin{equation}
\label{equ:derived_isomorphism}
  \projmod 
  \stackrel{\sim}{\longrightarrow} 
  \RHom(\injcomplx,\injcomplx \LTensor \projmod).
\end{equation}
Here 
\[
  \injcomplx \LTensor \projmod \cong \injcomplx \otimes \projmod
\]
because $\projmod$ is projective.  Moreover, $\injcomplx \otimes
\projmod$ is a bounded complex of injective modules so
\[
  \RHom(\injcomplx,\injcomplx \LTensor \projmod) 
  \cong \RHom(\injcomplx,\injcomplx \otimes \projmod) 
  \cong \Hom(\injcomplx,\injcomplx \otimes \projmod).  
\]
So the isomorphism \eqref{equ:derived_isomorphism} in the derived
category is represented by the chain map
\[
  \projmod 
  \longrightarrow
  \Hom(\injcomplx,\injcomplx \otimes \projmod)
\]
which must accordingly be a quasi-isomorphism.  

Completing to a distinguished triangle in $\sK(\generalR)$ gives
\[
  \projmod 
  \longrightarrow \Hom(\injcomplx,\injcomplx \otimes \projmod) 
  \longrightarrow C 
  \longrightarrow
\]
where $C$ is exact.  Here $\injcomplx$ and $\injcomplx \otimes
\projmod$ are bounded, so $\Hom(\injcomplx,\injcomplx \otimes
\projmod)$ is bounded.  As the same is true for $\projmod$, the
mapping cone $C$ is also bounded.

Now, the distinguished triangle gives another distinguished triangle
\[
  \Hom(\projcomplx,\projmod) 
  \longrightarrow
  \Hom(\projcomplx,\Hom(\injcomplx,\injcomplx \otimes \projmod)) 
  \longrightarrow
  \Hom(\projcomplx,C) 
  \longrightarrow.
\]
Here $\Hom(\projcomplx,C)$ is exact because $\projcomplx$ is a complex
of projective modules while $C$ is a bounded exact complex.  So to see
that $\Hom(\projcomplx,\projmod)$ is exact as desired, it is enough to
see that $\Hom(\projcomplx,\Hom(\injcomplx,\injcomplx \otimes
\projmod))$ is exact.

However,
\[
  \Hom(\projcomplx,\Hom(\injcomplx,\injcomplx \otimes \projmod)) 
  \cong 
  \Hom(\injcomplx \otimes \projcomplx,\injcomplx \otimes \projmod).
\]
And this is exact because $\injcomplx \otimes \projcomplx$ is exact by
assumption while $\injcomplx \otimes \projmod$ is a bounded complex of
injective modules.
\end{proof}

\begin{Lemma}
\label{lem:cg}
The homotopy category of complexes of projective \generalRm s,
$\KProjgeneralR$, is a compactly generated triangulated category.
\end{Lemma}

\begin{proof}
It is clear that $\KProjgeneralR$ is triangulated.

The ring $\generalR$ is noetherian and hence coherent, and by remark
\ref{rmk:flats_have_finite_pd} each flat \generalRm\ has finite projective
dimension.  So $\KProjgeneralR$ is compactly generated by 
\cite[thm.\ 2.4]{PJhomproj}.  
\end{proof}

Combining lemmas \ref{lem:technical} and \ref{lem:cg} with Bousfield
localization now gives existence of the adjoint functor $\adj$.

\begin{Proposition}
\label{pro:adjoint}
The inclusion functor $\inc : \EnochsgeneralR \longrightarrow
\KProjgeneralR$ has a right-adjoint $\adj : \KProjgeneralR
\longrightarrow \EnochsgeneralR$. 
\end{Proposition}

\begin{proof}
Consider the functor
\[
  k(-) = \H^0((\generalR \oplus \injcomplx) \otimes_{\generalR} -)
\]
from $\KProjgeneralR$ to $\Ab$, the category of abelian groups.  This
is clearly a homological functor respecting set indexed coproducts.
Moreover,
\[
  k(\Sigma^i \projcomplx) 
  \cong 
  \H^i(\projcomplx) \oplus \H^i(\injcomplx \otimes \projcomplx),
\]
where $\Sigma^i$ denotes $i$'th suspension, so for $\projcomplx$ to
satisfy $k(\Sigma^i \projcomplx) = 0$ for each $i$ means
\[
  \H^i(\projcomplx) = 0
\]
and
\[
  \H^i(\injcomplx \otimes \projcomplx) = 0
\]
for each $i$.  Using lemma \ref{lem:technical}, this shows
\[
  \{\, \projcomplx \in \KProjgeneralR
       \,\mid\, k(\Sigma^i \projcomplx) = 0 \mbox{ for each } i \,\}
  = \EnochsgeneralR.
\]
That is, $\EnochsgeneralR$ is the kernel of the homological functor $k$.

One consequence of this is that $\EnochsgeneralR$ is closed under set
indexed coproducts.  Hence \cite[lem.\ 3.5]{KrauseSmashing} says that
for $\inc$ to have a right-adjoint is the same as for the Verdier
quotient $\KProjgeneralR/\EnochsgeneralR$ to satisfy that each $\Hom$
set is in fact a set (as opposed to a class).

Now, the category $\KProjgeneralR$ is compactly generated by lemma
\ref{lem:cg}.  By \cite[lem.\ 4.5.13]{NeemanBook} with $\beta =
\aleph_0$, this even implies that there is only a set of isomorphism
classes of compact objects in $\KProjgeneralR$.  Hence the version of
Bousfield localization given in \cite[thm.\ 4.1]{PJSpectra} applies to
the functor $k$ on $\KProjgeneralR$, and gives that $\KProjgeneralR$
modulo the kernel of $k$ satisfies that each $\Hom$ is a set.  That
is, $\KProjgeneralR/\EnochsgeneralR$ satisfies that each $\Hom$ is a
set, as desired.
\end{proof}

The methods given above also apply to non-commutative algebras.  Let
the following setups replace setups \ref{set:blanket1} and
\ref{set:I1}.

\medskip
\noindent
{\bf Setup \ref{set:blanket1}'. } 
Let $\generalR$ be a left-coherent and right-noetherian $k$-algebra
over the field $k$ for which there exists a left-noetherian $k$-algebra
$\matlisgeneralR$ and a dualizing complex
${}_{\matlisgeneralR}D_{\generalR}$.  That is, $D$ is a complex of
$\matlisgeneralR$-left-$\generalR$-right-modules, and
\begin{enumerate}

  \item  The cohomology of $D$ is bounded and finitely generated
         both over $\matlisgeneralR$ and over $\generalR^{\opp}$.

  \item  The injective dimensions $\id_{\matlisgeneralR} D$ and
         $\id_{\generalR^{\opp}} D$ are finite.

  \item  The canonical morphisms 
\[
  \generalR \longrightarrow \RHom_{\matlisgeneralR}(D,D) 
  \;\;\; \mbox{and} \;\;\;
  \matlisgeneralR \longrightarrow \RHom_{\generalR^{\opp}}(D,D) 
\]
in the derived categories $\D(\generalR \otimes_k \generalR^{\opp})$
and $\D(\matlisgeneralR \otimes_k \matlisgeneralR^{\opp})$ are
isomorphisms.

\end{enumerate}

\medskip
\noindent
{\bf Setup \ref{set:I1}'. } 
Let $D \stackrel{\simeq}{\longrightarrow} \injcomplx$ be an injective
resolution of $D$ over $\matlisgeneralR \otimes_k \generalR^{\opp}$.
Replace $\injcomplx$ by a bounded truncation.  This may ruin the
property that $\injcomplx$ is an injective resolution over
$\matlisgeneralR \otimes_k \generalR^{\opp}$, but because
$\id_{\matlisgeneralR} D$ and $\id_{\generalR^{\opp}} D$ are finite, I
can still suppose that $\injcomplx$ consists of modules which are
injective both over $\matlisgeneralR$ and over $\generalR^{\opp}$.
\medskip

The definition of dualizing complexes over non-commutative algebras is
due to \cite[def.\ 1.1]{YekutieliZhang}.

With setups \ref{set:blanket1} and \ref{set:I1} replaced by setups
\ref{set:blanket1}' and \ref{set:I1}', let me inspect the rest of this
section.  As the ground ring $\generalR$ is now non-commutative, I
must replace ``module'' by ``left-module'' throughout.  Remark
\ref{rmk:flats_have_finite_pd} also needs to be replaced by the
following.

\medskip
\noindent
{\bf Remark \ref{rmk:flats_have_finite_pd}'. } 
Under setup \ref{set:blanket1}', each flat \generalRlm\ has finite
projective dimension by \cite{PJfdpd}.
\medskip

After this, the proof of lemma \ref{lem:technical} goes through if one
keeps track of left and right structures throughout, and the proofs of
lemma \ref{lem:cg} and proposition \ref{pro:adjoint} also still work.

Hence I can sum up the results of this section in the following
theorem. 

\begin{Theorem}
\label{thm:adjoint}
Consider either of the following two situations.
\begin{enumerate}

  \item  $\generalR$ is a noetherian commutative ring with a dualizing
         complex (see setup \ref{set:blanket1}).

  \medskip

  \item  $\generalR$ is a left-coherent and right-noetherian
         $k$-algebra over the field $k$ for which there exists a
         left-noetherian $k$-algebra $\matlisgeneralR$ and a dualizing
         complex ${}_{\matlisgeneralR}D_{\generalR}$ (see setup
         \ref{set:blanket1}'). 

\end{enumerate}
Then the inclusion functor 
\[
  \inc : \EnochsgeneralR \longrightarrow \KProjgeneralR
\]
has a right-adjoint 
\[
  \adj : \KProjgeneralR \longrightarrow \EnochsgeneralR.
\]
\end{Theorem}

\section{Gorenstein projective resolutions}
\label{sec:resolution}

This section shows the existence of Gorenstein projective resolutions
when the adjoint functor $\adj$ exists.

\begin{Setup}
\label{set:blanket}
For the rest of this paper, $\R$ is a ring for which the inclusion
functor
\[
  \inc : \EnochsR \longrightarrow \KProjR
\]
has a right-adjoint
\[
  \adj : \KProjR \longrightarrow \EnochsR.
\]
\end{Setup}

\begin{Remark}
\label{rmk:blanket}
The existence of the right-adjoint $\adj$ is precisely the hypothesis
under which the constructions of this paper work.  

The functor $\adj$ exists over fairly general rings; see theorem
\ref{thm:adjoint} and remark \ref{rmk:dualizing}.  As mentioned in the
introduction, it may even be the case that the functor $\adj$ exists
over {\em any} ring, but I do not know how to prove that.
\end{Remark}

\begin{Remark}
\label{rmk:approximation}
If $\projcomplx$ is a complex of projective \Rlm s, then $\adj
\projcomplx$ can be thought of as the best approximation to
$\projcomplx$ by a complete projective resolution.

Elaborating on this, if $\M$ is an \Rlm\ with projective
resolution $\projcomplx$, then $\adj \projcomplx$ can be thought of as
the best approximation to $\M$ by a complete projective
resolution.  This point will be made more precise in lemma
\ref{lem:adj_gives_complete_resolutions}.
\end{Remark}

\begin{Construction}
\label{con:surjective}
If $\projcomplx$ is a complex of \Rlm s, then for each $i$ there is a
chain map
\[
  \begin{diagram}[labelstyle=\scriptstyle,width=4ex,height=7ex]
    \cdots & \rTo & 0 & & \rTo & & \projcomplx^i & & \rTo^{\id} & & \projcomplx^i & & \rTo & & 0 & \rTo & \cdots \\
    & & \dTo & & & & \dTo^{\id} & & & & \dTo_{\partial^i} & & & & \dTo & & \\
    \cdots & \rTo & \projcomplx^{i-1} & & \rTo_{\partial^{i-1}} & & \projcomplx^i & & \rTo_{\partial^i} & & \projcomplx^{i+1} & & \rTo_{\partial^{i+1}} & & \projcomplx^{i+2} & \rTo & \cdots \\
  \end{diagram}
\]
where the upper complex is null homotopic.  

If $\Tatecomplx \stackrel{\Tatecomplxmap}{\longrightarrow}
\projcomplx$ is now a chain map, then I can add the upper complex to
$\Tatecomplx$ and thereby change $\Tatecomplxmap$ so that the $i$'th
component $\Tatecomplx^i \stackrel{\Tatecomplxmap^i}{\longrightarrow}
\projcomplx^i$ becomes surjective.  Doing so does not change the
isomorphism class of $\Tatecomplxmap$ in $\sK(\R)$, the homotopy
category of complexes of \Rlm s.
\end{Construction}

\begin{Construction}
\label{con:basic_exact_sequence}
If $\M$ is an \Rlm, then let $\projcomplx$ be a projective resolution
concentrated in non-positive cohomological degrees and consider the
counit morphism $\inc \adj \projcomplx
\stackrel{\epsilon_{\projcomplx}}{\longrightarrow} \projcomplx$ in $\KProjR$.
By applying construction \ref{con:surjective} in each degree, I can
assume that $\epsilon_{\projcomplx}$ is represented by a surjective
chain map, so setting
\[
  \Enochscomplxtil = \inc \adj \projcomplx, 
\]
there is a short exact sequence of complexes
\[
  0 \rightarrow \kercomplx 
    \longrightarrow \Enochscomplxtil
    \longrightarrow \projcomplx \rightarrow 0.
\]
Note that since both $\Enochscomplxtil$ and $\projcomplx$ consist of
projective modules, the sequence is semi-split (that is, split in each
degree) and $\kercomplx$ also consists of projective modules.
\end{Construction}

\begin{Lemma}
\label{lem:zero}
Consider the complex $\kercomplx$ from construction
\ref{con:basic_exact_sequence}.  Then
\[
  \Hom_{\KProjRsmall}(\Enochscomplx,\kercomplx) = 0
\]
for $\Enochscomplx$ in $\EnochsR$.
\end{Lemma}

\begin{proof}
The chain map $\Enochscomplxtil \longrightarrow \projcomplx$
represents the counit morphism
\[
  \inc \adj \projcomplx
  \stackrel{\epsilon_{\projcomplx}}{\longrightarrow} 
  \projcomplx
\]
which leads to a commutative diagram
\[
  \begin{diagram}[labelstyle=\scriptstyle,width=12ex,height=5ex]
    \Hom_{\EnochsRsmall}(\Enochscomplx,\adj \projcomplx)
      & \rTo^{\inc(-)}
      & \Hom_{\KProjRsmall}(\inc \Enochscomplx,\inc \adj \projcomplx) \\
    & \SE & \dTo_{\Hom(\inc \Enochscomplx,\epsilon_{\projcomplx})} \\
    & & \Hom_{\KProjRsmall}(\inc \Enochscomplx,\projcomplx) \lefteqn{,} \\
  \end{diagram}
\]
where the diagonal map is the adjunction isomorphism while the
horizontal map is an isomorphism because $\inc$ is the inclusion
functor of a full subcategory.  The vertical map must therefore also
be an isomorphism.  That is,
\begin{equation}
\label{equ:FPiso}
  \Hom_{\KProjRsmall}(\Enochscomplx,\Enochscomplxtil)
  \longrightarrow \Hom_{\KProjRsmall}(\Enochscomplx,\projcomplx)
\end{equation}
is an isomorphism.

Now, the short exact sequence from construction
\ref{con:basic_exact_sequence} is semi-split and therefore gives a
distinguished triangle
\[
  \kercomplx 
  \longrightarrow \Enochscomplxtil
  \longrightarrow \projcomplx 
  \longrightarrow
\]
in $\KProjR$.  Hence there is a long exact sequence consisting of
pieces
\[
  {\scriptstyle 
  \Hom_{\KProjRsmall}(\Sigma^i \Enochscomplx,\kercomplx)
  \longrightarrow 
  \Hom_{\KProjRsmall}(\Sigma^i \Enochscomplx,\Enochscomplxtil)
  \longrightarrow
  \Hom_{\KProjRsmall}(\Sigma^i \Enochscomplx,\projcomplx).
  }
\]
Since $\Sigma^i \Enochscomplx$ is in $\EnochsR$ for each $i$, the
second homomorphism here is of the type from equation
\eqref{equ:FPiso}, so is an isomorphism for each $i$.  This implies
$\Hom_{\KProjRsmall}(\Enochscomplx,\kercomplx) = 0$ as desired.
\end{proof}

\begin{Remark}
\label{rmk:Gorproj}
For the following lemma, recall that a Gorenstein projective \Rlm\ is
a module of the form $\Gorproj = \Ker(\Enochscomplx^1 \longrightarrow
\Enochscomplx^2)$ where $\Enochscomplx$ is in $\EnochsR$ (cf.\
definition \ref{def:Enochs}).
\end{Remark}

\begin{Lemma}
\label{lem:Gorenstein_projective_epic}
Consider the complex $\kercomplx$ from construction
\ref{con:basic_exact_sequence}.  Suppose that the sequence
\[
  \cdots \longrightarrow \kercomplx^{i-2} \longrightarrow \kercomplx^{i-1}
  \stackrel{k}{\longrightarrow} \N \rightarrow 0
\]
obtained from $\kercomplx$ is exact.  

Let $\Gorprojtil$ be Gorenstein projective and let $\Gorprojtil
\stackrel{\Gorprojtilmap}{\longrightarrow} \N$ be a homomorphism.
Then $\Gorprojtilmap$ lifts through $k$,
\[
  \begin{diagram}[labelstyle=\scriptstyle]
                     &           & \Gorprojtil \\
                     & \ldDotsto & \dTo_{\Gorprojtilmap} \\
    \kercomplx^{i-1} & \rTo_{k}  & \N \lefteqn{.}
  \end{diagram}
\]
\end{Lemma}

\begin{proof}
By (de)suspending, I can clearly pick a complex $\Enochscomplx$ in
$\EnochsR$ with $\Gorprojtil = \Ker(\Enochscomplx^i \longrightarrow
\Enochscomplx^{i+1})$, and it is not hard to see that there is a chain
map $\Enochscomplx \stackrel{\Enochscomplxmap}{\longrightarrow}
\kercomplx$ which fits together with $\Gorprojtil
\stackrel{\Gorprojtilmap}{\longrightarrow} \N$ in a commutative
diagram 
\[
  \begin{diagram}[labelstyle=\scriptstyle,height=4ex,width=4ex]
    \cdots & \rTo & \Enochscomplx^{i-2} & & \rTo & & \Enochscomplx^{i-1} & & \rTo & & \Enochscomplx^i & & \rTo & & \Enochscomplx^{i+1} & \rTo & \cdots \\
    & & & & & & & \rdOnto & & \ruEmbed_{\ell} & & & & & & & \\
    & & \dTo^{\Enochscomplxmap^{i-2}} & & & & \dTo^{\Enochscomplxmap^{i-1}} & & \Gorprojtil & & \dTo_{\Enochscomplxmap^i} & & & & \dTo_{\Enochscomplxmap^{i+1}} & & \\
    & & & & & & & & \vLine^{\Gorprojtilmap} & & & & & & & & \\
    \cdots & \rTo & \kercomplx^{i-2} & & \rTo & & \kercomplx^{i-1} & \rTo & \HonV & & \kercomplx^i & & \rTo & & \kercomplx^{i+1} & \rTo & \cdots \lefteqn{.} \\
    & & & & & & & \rdOnto_k & \dTo & \ruTo & & & & & & & \\
    & & & & & & & & \N & & & & & & & & \\
  \end{diagram}
\]

Since lemma \ref{lem:zero} says
$\Hom_{\KProjRsmall}(\Enochscomplx,\kercomplx) = 0$ for
$\Enochscomplx$ in $\EnochsR$, the chain map $\Enochscomplxmap$ must
be null homotopic.  Let $\epsilon$ be a null homotopy with $\Enochscomplxmap =
\epsilon\partial^{\Enochscomplx} + \partial^{\kercomplx}\epsilon$,
consisting of components $\Enochscomplx^j
\stackrel{\epsilon^j}{\longrightarrow}
\kercomplx^{j-1}$.  Then it is straightforward to prove
\[
  k \circ (\epsilon^i \ell) = \Gorprojtilmap,
\]
so $\Gorprojtil \stackrel{\Gorprojtilmap}{\longrightarrow} \N$
has been lifted through $\kercomplx^{i-1} \stackrel{k}{\longrightarrow}
\N$ as desired.
\end{proof}

\begin{Remark}
\label{rmk:Gorenstein_projective_resolutions}
For the next theorem, recall that an augmented Gorenstein projective
resolution of an \Rlm\ $\M$ is an exact sequence
\[
  \cdots \longrightarrow \Gorproj_2
  \longrightarrow \Gorproj_1
  \longrightarrow \Gorproj_0 
  \longrightarrow \M
  \rightarrow 0
\]
where the $\Gorproj_i$ are Gorenstein projective modules,
which stays exact when one applies the functor
$\Hom(\Gorprojtil,-)$ for any Gorenstein projective module
$\Gorprojtil$.  The complex
\[
  \Gorproj = 
  \cdots \longrightarrow \Gorproj_2
  \longrightarrow \Gorproj_1
  \longrightarrow \Gorproj_0
  \longrightarrow 0
  \longrightarrow \cdots
\]
is then called a Gorenstein projective resolution of $\M$.  
\end{Remark}

\begin{Remark}
\label{rmk:resolution}
Recall construction \ref{con:basic_exact_sequence}.  The complex
$\Enochscomplxtil$ is in $\EnochsR$.  In particular it is exact, and
therefore the cohomology long exact sequence shows
\[
  \H^i\!\kercomplx =
  \left\{
    \begin{array}{cl}
      \M & \mbox{ for } i =     1, \\
      0       & \mbox{ for } i \not= 1.
    \end{array}
  \right.
\]
Hence there is an exact sequence
\[
  \cdots 
  \longrightarrow \kercomplx^{-2}
  \longrightarrow \kercomplx^{-1}
  \longrightarrow \kercomplx^0
  \longrightarrow \Ker \partial_{\kercomplx}^1
  \longrightarrow \M
  \rightarrow 0.
\]
\end{Remark}

\begin{Theorem}
\label{thm:resolution}
Let $\M$ be an \Rlm.  Then the exact sequence from remark
\ref{rmk:resolution} is an augmented Gorenstein projective resolution
of $\M$.
\end{Theorem}

\begin{proof}
The modules $\kercomplx^0, \kercomplx^{-1}, \ldots$ are projective and
hence Gorenstein projective.

As for $\Ker \partial_{\kercomplx}^1$, observe that in the short exact
sequence from construction \ref{con:basic_exact_sequence}, the complex
$\projcomplx$ is concentrated in non-positive cohomological degrees,
so the modules $P^1$ and $P^2$ are zero.  So in degrees $1$ and $2$,
the short exact sequence gives
\[
  \begin{diagram}[labelstyle=\scriptstyle,midshaft]
    \kercomplx^2 & \rTo^{\cong} & \Enochscomplxtil^2 \\
    \uTo^{\partial_{\kercomplx}^1} & & \uTo_{\partial_{\Enochscomplxtil}^1} \\
    \kercomplx^1 & \rTo^{\cong} & \Enochscomplxtil^1 \lefteqn{.} \\
  \end{diagram}
\]
Hence $\Ker \partial_{\kercomplx}^1 \cong \Ker
\partial_{\Enochscomplxtil}^1$, and $\Ker \partial_{\Enochscomplxtil}^1$
is Gorenstein projective because $\Enochscomplxtil$ is in
$\EnochsR$.

To complete the proof, I must show that the exact sequence from remark
\ref{rmk:resolution},
\[
  \cdots 
  \longrightarrow \kercomplx^{-2}
  \longrightarrow \kercomplx^{-1}
  \longrightarrow \kercomplx^0
  \longrightarrow \Ker \partial_{\kercomplx}^1
  \longrightarrow \M
  \rightarrow 0,
\]
remains exact when one applies the functor
$\Hom(\Gorprojtil,-)$ for any Gorenstein projective module
$\Gorprojtil$.

First, let $i \leq 0$ be an integer and let $\Gorprojtil
\stackrel{\Gorprojtilmap}{\longrightarrow} \kercomplx^i$ be a
homomorphism whose composition with the subsequent homomorphism in the
exact sequence is zero.  I must show that $\Gorprojtilmap$ lifts through
$\kercomplx^{i-1} \longrightarrow \kercomplx^i$.

I can view $\Gorprojtilmap$ as a homomorphism $\Gorprojtil
\stackrel{\Gorprojtilmap}{\longrightarrow} \Ker
\partial_{\kercomplx}^i$, and must then show that $\Gorprojtilmap$
lifts through the canonical homomorphism $\kercomplx^{i-1}
\longrightarrow \Ker \partial_{\kercomplx}^i$.  But this
follows from lemma \ref{lem:Gorenstein_projective_epic} applied to
\[
  \cdots \longrightarrow \kercomplx^{i-2}
  \longrightarrow \kercomplx^{i-1}
  \longrightarrow \Ker \partial_{\kercomplx}^i
  \rightarrow 0.
\]

Secondly, let $\Gorprojtil \stackrel{\Gorprojtilmap}{\longrightarrow}
\Ker \partial_{\kercomplx}^1$ be a homomorphism whose composition with
the subsequent homomorphism in the exact sequence, $\Ker
\partial_{\kercomplx}^1 \longrightarrow \M$, is zero.  I must show
that $\Gorprojtilmap$ lifts through $\kercomplx^0 \longrightarrow \Ker
\partial_{\kercomplx}^1$.  

I can view $\Gorprojtilmap$ as a homomorphism $\Gorprojtil
\stackrel{\Gorprojtilmap}{\longrightarrow} \Image
\partial_{\kercomplx}^0$, and must then show that $\Gorprojtilmap$
lifts through the canonical homomorphism $\kercomplx^0 \longrightarrow
\Image \partial_{\kercomplx}^0$.  But this follows from lemma
\ref{lem:Gorenstein_projective_epic} applied to 
\[
  \cdots \longrightarrow \kercomplx^{-1}
  \longrightarrow \kercomplx^0
  \longrightarrow \Image \partial_{\kercomplx}^0
  \rightarrow 0.
\]

Thirdly, let $\Gorprojtil
\stackrel{\Gorprojtilmap}{\longrightarrow} \M$ be a homomorphism.
I must show that $\Gorprojtilmap$ lifts through $\Ker
\partial_{\kercomplx}^1 \longrightarrow \M$.
However, from the data given I can construct a commutative diagram
\[
  \begin{diagram}[labelstyle=\scriptstyle,height=4ex,width=4ex]
    \cdots & \rTo & \kercomplx^0 & & \rTo^{\partial_{\kercomplx}^0} & & \kercomplx^1 & & \rTo & & \kercomplx^2 & \rTo & \cdots \lefteqn{,} \\
    & & & \rdTo & & \ruEmbed & & \rdOnto & & \ruTo & & & \\
    & & & & \Ker \partial_{\kercomplx}^1 & & & & \Coker \partial_{\kercomplx}^0 & & & & \\
    & & & & & \rdOnto & & \ruEmbed_j & & & & & \\
    & & & & & & \M & & & & & & \\
  \end{diagram}
\]
and by applying lemma \ref{lem:Gorenstein_projective_epic} to
\[
  \cdots \longrightarrow \kercomplx^0
  \longrightarrow \kercomplx^1
  \longrightarrow \Coker \partial_{\kercomplx}^0
  \rightarrow 0
\]
I find that $\Gorprojtil \stackrel{j \Gorprojtilmap}{\longrightarrow}
\Coker \partial_{\kercomplx}^0$ lifts through $\kercomplx^1
\longrightarrow \Coker \partial_{\kercomplx}^0$.  It is a small
diagram exercise to see that hence, $\Gorprojtil
\stackrel{\Gorprojtilmap}{\longrightarrow} \M$ lifts through $\Ker
\partial_{\kercomplx}^1 \longrightarrow \M$ as desired.
\end{proof}

Let me close the section with the following easy consequence.

\begin{Remark} 
Recall that a Gorenstein projective precover of an \Rlm\ $\M$ is a
homomorphism $\Gorproj \stackrel{\Gorprojmap}{\longrightarrow} \M$
where $\Gorproj$ is a Gorenstein projective module, so that if
$\Gorprojtil$ is any Gorenstein projective module with a homomorphism
$\Gorprojtil \stackrel{\Gorprojtilmap}{\longrightarrow} \M$, then
$\Gorprojtilmap$ lifts through $\Gorprojmap$,
\[
  \begin{diagram}[labelstyle=\scriptstyle]
             &                     & \Gorprojtil \\
             & \ldDotsto           & \dTo_{\Gorprojtilmap} \\
    \Gorproj & \rTo_{\Gorprojmap}  & \M \lefteqn{.}
  \end{diagram}
\smallskip
\]
\end{Remark}

\begin{Theorem}
Each \Rlm\ has a Gorenstein projective precover.
\end{Theorem}

\begin{proof}
It follows from theorem \ref{thm:resolution} that the homomorphism
\[
  \Ker \partial_{\kercomplx}^1 \longrightarrow \M 
\]
is a Gorenstein projective precover.
\end{proof}

\section{Tate $\Ext$ groups}
\label{sec:Tate}

This section defines Tate $\Ext$ groups, and goes on to show some simple
properties: A short exact sequence in either variable gives rise to a
long exact sequence of Tate $\Ext$ groups; when the Tate $\Ext$ groups from
\cite{AvrMart} and \cite{Veliche} are defined, they agree with the
ones defined in this paper; and classical Tate cohomology is the
special case $\TExt_{kG}^i(k,\N)$ of the Tate $\Ext$ groups.  Finally, it
is proved that the Tate $\Ext$ groups fit into the long exact sequence
\eqref{equ:long_exact_sequence} from the introduction.

\begin{Remark}
It is classical that the category of \Rlm s $\Mod(\R)$ is equivalent
to the full subcategory of $\KProjR$ consisting of projective
resolutions of \Rlm s.  Let 
\[
  \res : \Mod(\R) \longrightarrow \KProjR
\]
be a functor implementing the equivalence.
\end{Remark}

\begin{Definition}
\label{def:Tate}
If $\M$ and $\N$ are \Rlm s, then the $i$'th Tate $\Ext$ group is
\[
  \TExt^i(\M,\N) = \H^i\!\Hom_{\R}(\adj \res \M,\N).
\]
\end{Definition}

\begin{Remark}
As pointed out in remark \ref{rmk:approximation}, the complex $\adj
\res \M$ can be thought of as the best approximation to $\M$ by a
complete projective resolution.  So taking $\Hom$ into $\N$ and taking
cohomology is the obvious way to get Tate $\Ext$ groups.
\end{Remark}

\begin{Proposition}
\label{pro:long_exact_sequences}
Let
\[
  0 \rightarrow \M^{\prime} 
    \longrightarrow \M 
    \longrightarrow \M^{\prime \prime} \rightarrow 0 
\]
and 
\[
  0 \rightarrow \N^{\prime}
    \longrightarrow \N 
    \longrightarrow \N^{\prime \prime} \rightarrow 0
\]
be short exact sequences of \Rlm s.  Then there are natural long exact
sequences
\[
  \cdots \longrightarrow \TExt^i(\M^{\prime \prime},\N)
  \longrightarrow \TExt^i(\M,\N)
  \longrightarrow \TExt^i(\M^{\prime},\N)
  \longrightarrow \cdots
\]
and
\[
  \cdots \longrightarrow \TExt^i(\M,\N^{\prime})
  \longrightarrow \TExt^i(\M,\N)
  \longrightarrow \TExt^i(\M,\N^{\prime \prime})
  \longrightarrow \cdots.
\]
\end{Proposition}

\begin{proof}
It is well known that the first short exact sequence in the
proposition results in a distinguished triangle in $\KProjR$,
\[
  \res \M^{\prime} 
  \longrightarrow \res \M
  \longrightarrow \res \M^{\prime \prime}
  \longrightarrow.
\]
Since $\inc$ is a triangulated functor, so is its adjoint $\adj$, so
there is also a distinguished triangle in $\EnochsR$,
\[
  \adj \res \M^{\prime} 
  \longrightarrow \adj \res \M
  \longrightarrow \adj \res \M^{\prime \prime}
  \longrightarrow.
\]
This again results in a distinguished triangle
\[
  {\scriptstyle
  \Hom_{\R}(\adj \res \M^{\prime \prime},\N)
  \longrightarrow \Hom_{\R}(\adj \res \M,\N)
  \longrightarrow \Hom_{\R}(\adj \res \M^{\prime},\N)
  \longrightarrow
  }
\]
whose cohomology long exact sequence is the first long exact sequence
in the proposition. 

The complex $\adj \res \M$ is in $\EnochsR$ so consists of projective
modules, so the second short exact sequence in the proposition gives a
short exact sequence of complexes
\[
  {\scriptstyle
  0 \rightarrow \Hom_{\R}(\adj \res \M,\N^{\prime})
    \longrightarrow \Hom_{\R}(\adj \res \M,\N)
    \longrightarrow \Hom_{\R}(\adj \res \M,\N^{\prime \prime})
    \rightarrow 0
  }
\]
whose cohomology long exact sequence is the second long exact sequence
in the proposition. 
\end{proof}

\begin{Remark}
\label{rmk:AvrMart_Tate}
If $\M$ and $\N$ are \Rlm s, then the earlier definition of Tate
$\Ext$ groups given in \cite{AvrMart} and \cite{Veliche} is
\[
  \TExt^i(\M,\N) = \H^i\!\Hom_{\R}(\Tatecomplx,\N)
\]
where $\Tatecomplx$ is a complete projective resolution of $\M$.  This
means that $\Tatecomplx$ is in $\EnochsR$ and sits in a diagram of
chain maps
\begin{equation}
\label{equ:complete_resolution}
  \Tatecomplx \stackrel{\Tatecomplxmap}{\longrightarrow} 
  \projcomplx \longrightarrow \M
\end{equation}
where $\projcomplx \longrightarrow \M$ is a projective resolution and
where $\Tatecomplx^i \stackrel{\Tatecomplxmap^i}{\longrightarrow}
\projcomplx^i$ is bijective for $i \ll 0$.

Note that not all \Rlm s have complete projective resolutions.  In
fact, the ones that do are exactly the ones which have finite
Gorenstein projective dimension by \cite[thm.\ 3.4]{Veliche}.
\end{Remark}

\begin{Lemma}
\label{lem:adj_gives_complete_resolutions}
Let $\M$ be an \Rlm\ which has a projective resolution $\projcomplx$
and a complete projective resolution $\Tatecomplx$.  Then
\[
  \adj \projcomplx \cong \Tatecomplx  
\]
in $\KProjR$.
\end{Lemma}

\begin{proof} 
All projective resolutions of $\M$ are isomorphic in $\KProjR$,
so I may as well prove the lemma for the specific projective
resolution $\projcomplx$ from equation
\eqref{equ:complete_resolution}.

By applying construction \ref{con:surjective} to the chain map
$\Tatecomplx \stackrel{\Tatecomplxmap}{\longrightarrow}
\projcomplx$ in cohomological degrees larger than some number, I can
assume that $\Tatecomplxmap$ is surjective.  Hence there is a short
exact sequence of complexes
\begin{equation}
\label{equ:pre_basic_exact_sequence}
  0 \rightarrow \kercomplx 
    \longrightarrow \Tatecomplx
    \stackrel{\Tatecomplxmap}{\longrightarrow} \projcomplx
    \rightarrow 0.
\end{equation}

Since both $\Tatecomplx$ and $\projcomplx$ consist of projective
modules, the sequence is semi-split and $\kercomplx$ also consists of
projective modules.  Moreover, by assumption, 
$\Tatecomplx^i \stackrel{\Tatecomplxmap^i}{\longrightarrow} \projcomplx^i$ 
is bijective for $i \ll 0$, so $K^i = 0$ for $i \ll 0$.  So
$\kercomplx$ is a left-bounded complex of projective modules.

Now let $\Enochscomplx$ be in $\EnochsR$.  In particular,
$\Hom_{\R}(\Enochscomplx,\projmod)$ is exact when $\projmod$ is a
projective module.  It is classical that
$\Hom_{\R}(\Enochscomplx,\kercomplx)$ is then also exact, because
$\kercomplx$ is a left-bounded complex of projective modules.  Indeed,
this follows by an argument analogous to the one which shows that if
$X$ is an exact complex and $I$ is a left-bounded complex of injective
modules, then $\Hom_{\R}(X,I)$ is exact.

Since the sequence \eqref{equ:pre_basic_exact_sequence} is semi-split,
it stays exact under the functor $\Hom_{\R}(\Enochscomplx,-)$.  So
there is a short exact sequence of complexes
\[
  0 \rightarrow \Hom_{\R}(\Enochscomplx,\kercomplx)
    \longrightarrow \Hom_{\R}(\Enochscomplx,\Tatecomplx)
    \longrightarrow \Hom_{\R}(\Enochscomplx,\projcomplx)
    \rightarrow 0.
\]
Since $\Hom_{\R}(\Enochscomplx,\kercomplx)$ is exact, the
cohomology long exact sequence shows that there is an isomorphism
\[
  \H^0\!\Hom_{\R}(\Enochscomplx,\Tatecomplx)
  \cong \H^0\!\Hom_{\R}(\Enochscomplx,\projcomplx)
\]
which is natural in $\Enochscomplx$.  That is, there is a natural
isomorphism
\[
  \Hom_{\KProjRsmall}(\Enochscomplx,\Tatecomplx)
  \cong \Hom_{\KProjRsmall}(\Enochscomplx,\projcomplx)
\]
which can also be written 
\[
  \Hom_{\EnochsRsmall}(\Enochscomplx,\Tatecomplx)
  \cong \Hom_{\KProjRsmall}(\Enochscomplx,\projcomplx)
\]
because $\Enochscomplx$ and $\Tatecomplx$ are in $\EnochsR$.

On the other hand, I also have a natural isomorphism
\[
  \Hom_{\KProjRsmall}(\Enochscomplx,\projcomplx)
  = \Hom_{\KProjRsmall}(\inc \Enochscomplx,\projcomplx)
  \cong \Hom_{\EnochsRsmall}(\Enochscomplx,\adj \projcomplx).
\]
Combining the last two equations gives a natural isomorphism
\[
  \Hom_{\EnochsRsmall}(\Enochscomplx,\Tatecomplx)
  \cong \Hom_{\EnochsRsmall}(\Enochscomplx,\adj \projcomplx),
\]
proving $\Tatecomplx \cong \adj \projcomplx$ as desired.
\end{proof}

\begin{Proposition}
\label{pro:correspondence}
Let $\M$ be an \Rlm\ which has a complete projective resolution
$\Tatecomplx$.  Then the Tate $\Ext$ groups of this paper
(see definition \ref{def:Tate}) coincide with the Tate $\Ext$ groups
which were defined in \cite{AvrMart} and \cite{Veliche} (see remark
\ref{rmk:AvrMart_Tate}).
\end{Proposition}

\begin{proof}
Lemma \ref{lem:adj_gives_complete_resolutions} gives that the
projective resolution $\res \M$ of $\M$ satisfies $\adj \res \M \cong
\Tatecomplx$.  Combining this with the formulae in definition
\ref{def:Tate} and remark \ref{rmk:AvrMart_Tate} proves the
proposition.  
\end{proof}

\begin{Proposition}
\label{pro:classical_Tate}
Let $k$ be a field, $G$ a finite group, and $\N$ a finite
dimensional $k$-linear representation of $G$.  Then the Tate $\Ext$ group
\[
  \TExt_{kG}^i(k,\N)
\]
of this paper is defined and isomorphic to the $i$'th classical Tate
cohomology group of $N$.
\end{Proposition}

\begin{proof}
The group algebra $kG$ is a finite dimensional $k$-algebra.  It is
clearly left-coherent and right-noetherian, and since it is in fact
self injective, it is clear that ${}_{kG}kG_{kG}$ is a
dualizing complex (cf.\ setup \ref{set:blanket1}').

Hence $\adj$ exists over $kG$ by theorem \ref{thm:adjoint}, and so the
Tate $\Ext$ groups of this paper are defined over $kG$.

The Tate $\Ext$ groups $\TExt_{kG}^i(k,\N)$ from \cite{AvrMart} and
\cite{Veliche} are also defined, and the $i$'th one is isomorphic to
the $i$'th classical Tate cohomology group of $\N$ according to
\cite[exam.\ 5.1]{AvrMart} with $k$ in place of $\BZ$.

But the Tate $\Ext$ groups of this paper and the ones from
\cite{AvrMart} and \cite{Veliche} are isomorphic by proposition
\ref{pro:correspondence}, so the present proposition follows.
\end{proof}

\begin{Definition}
If $\M$ and $\N$ are \Rlm s, then the $i$'th Gorenstein $\Ext$ group
$\GExt^i(\M,\N)$ is
\[
  \GExt^i(\M,\N) = \H^i \Hom_{\R}(\Gorproj,\N)
\]
where $\Gorproj$ is a Gorenstein projective resolution of $\M$
(cf.\ remark \ref{rmk:Gorenstein_projective_resolutions}).  
\end{Definition}

\begin{Remark}
The resolution $\Gorproj$ exists by theorem \ref{thm:resolution}.
Note that $\GExt^i(-,-)$ is a well defined bifunctor; see
\cite{AvrMart} or \cite{HHGorensteinDerived} for this and other
properties.  
\end{Remark}

\begin{Construction}
\label{con:truncated_sequence}
Consider the short exact sequence from construction
\ref{con:basic_exact_sequence},
\[
  0 \rightarrow \kercomplx 
    \longrightarrow \Enochscomplxtil
    \longrightarrow \projcomplx \rightarrow 0,
\]
where $\projcomplx$ is a projective resolution of the \Rlm\ $\M$ and
where $\Enochscomplxtil = \inc \adj \projcomplx$.  Truncating the
complexes $\kercomplx$ and $\Enochscomplxtil$ gives a new short exact
sequence of complexes,
\[
  \begin{diagram}[labelstyle=\scriptstyle,height=4ex,width=7ex]
    & & \vdots & & \vdots & & \vdots & & \\
    & & \uTo & & \uTo & & \uTo & & \\
    0 & \rTo & 0 & \rTo & 0 & \rTo & 0 & \rTo & 0 \\
    & & \uTo & & \uTo & & \uTo & & \\
    0 & \rTo & \Ker \partial_{\kercomplx}^1 & \rTo & \Ker \partial_{\Enochscomplxtil}^1 & \rTo & 0 & \rTo & 0 \\
    & & \uTo & & \uTo & & \uTo & & \\
    0 & \rTo & \kercomplx^0 & \rTo & \Enochscomplxtil^0 & \rTo & \projcomplx^0 & \rTo & 0 \\
    & & \uTo & & \uTo & & \uTo & & \\
    0 & \rTo & \kercomplx^{-1} & \rTo & \Enochscomplxtil^{-1} & \rTo & \projcomplx^{-1} & \rTo & 0 \lefteqn{,} \\
    & & \uTo & & \uTo & & \uTo & & \\
    & & \vdots & & \vdots & & \vdots & & \\
  \end{diagram}
\]
which I will denote
\begin{equation}
\label{equ:truncated_sequence}
  0 \rightarrow \kercomplxprime \longrightarrow \Enochscomplxtilprime
  \longrightarrow \projcomplx \rightarrow 0.
\end{equation}
\end{Construction}

\begin{Theorem}
\label{thm:long_exact_sequence}
Let $\M$ and $\N$ be \Rlm s.  Then there is a long exact
sequence 
\begin{eqnarray*}
  0 & \longrightarrow & \GExt^1(\M,\N) 
                        \longrightarrow \Ext^1(\M,\N)
                        \longrightarrow \TExt^1(\M,\N) \\
    & \longrightarrow & \cdots \\
    & \longrightarrow & \GExt^i(\M,\N) 
                        \longrightarrow \Ext^i(\M,\N)
                        \longrightarrow \TExt^i(\M,\N) 
                        \longrightarrow \cdots, \\
\end{eqnarray*}
natural in $\M$ and $\N$.
\end{Theorem}

\begin{proof}
Consider the short exact sequence \eqref{equ:truncated_sequence} from
construction \ref{con:truncated_sequence}.  The complex $\projcomplx$
is a projective resolution of $\M$ and in order to make everything
natural in $\M$, I can clearly suppose
\[
  \projcomplx = \res \M
\]
where $\res \M$ is a projective resolution depending functorially on
$\M$.  Since $\projcomplx = \res \M$ consists of projective modules,
the short exact sequence \eqref{equ:truncated_sequence} is semi-split
and therefore stays exact under the functor $\Hom_{\R}(-,\N)$.  So
there is a short exact sequence of complexes
\begin{equation}
\label{equ:Hom_of_truncated_sequence}
  0 \rightarrow \Hom_{\R}(\res \M,\N)
    \longrightarrow \Hom_{\R}(\Enochscomplxtilprime,\N)
    \longrightarrow \Hom_{\R}(\kercomplxprime,\N)
    \rightarrow 0.
\end{equation}

Since $\res \M$ is a projective resolution of $\M$, I have
\[
  \H^i\!\Hom_{\R}(\res \M,\N) = \Ext^i(\M,\N)
\]
for each $i$.  

The complex 
\[
  \Enochscomplxtil = \inc \adj \projcomplx = \adj \res \M
\]
is in $\EnochsR$ so it is exact, so
\[
  \Enochscomplxtilprime =
  \cdots 
  \longrightarrow \Enochscomplxtil^{-1}
  \longrightarrow \Enochscomplxtil^0
  \longrightarrow \Ker \partial_{\Enochscomplxtil}^1
  \longrightarrow 0
  \longrightarrow \cdots
\]
is also exact, and hence $\H^0\!\Hom_{\R}(\Enochscomplxtilprime,\N) =
0$.  On the other hand, the form of $\Enochscomplxtilprime$ makes it
clear that
\begin{eqnarray*}
  \H^i\!\Hom_{\R}(\Enochscomplxtilprime,\N)
  & = & \H^i\!\Hom_{\R}(\Enochscomplxtil,\N) \\
  & = & \H^i\!\Hom_{\R}(\adj \res \M,\N) \\
  & = & \TExt^i(\M,\N)
\end{eqnarray*}
for $i \geq 1$.

Finally, theorem \ref{thm:resolution} says that 
\[
  \kercomplxprime =
  \cdots 
  \longrightarrow \kercomplx^{-1}
  \longrightarrow \kercomplx^0
  \longrightarrow \Ker \partial_{\kercomplx}^1
  \longrightarrow 0
  \longrightarrow \cdots
\]
is a Gorenstein projective resolution of $\M$, shifted one step to the
right.  Hence
\[
  \H^i\!\Hom_{\R}(\kercomplxprime,\N)
  = \GExt^{i+1}(\M,\N)
\]
for $i \geq -1$.  

So looking at the cohomology long exact sequence of
\eqref{equ:Hom_of_truncated_sequence}, starting with
$\H^0\!\Hom_{\R}(\Enochscomplxtilprime,\N) = 0$, gives
\[
  0 \rightarrow \GExt^1(\M,\N) 
    \longrightarrow \Ext^1(\M,\N)
    \longrightarrow \TExt^1(\M,\N)
    \longrightarrow \cdots
\]
as desired.
\end{proof}

\bigskip

\noindent
{\bf Acknowledgement.}  
I would like to thank Edgar Enochs and Amnon Neeman for a number of
illuminating comments and answers.

This paper supersedes two earlier manuscripts, ``The Gorenstein
projective modules are precovering'' and ``Tate cohomology over fairly
general rings''.

The diagrams were typeset with Paul Taylor's {\tt diagrams.tex}.

\end{document}